\documentclass[12pt,final]{l4dc2020} 

\title[RKHS embedding for stochastic programming and control]{A Kernel Mean Embedding Approach to Reducing Conservativeness\\
in Stochastic Programming and Control}
\usepackage{times}

\usepackage{amsmath}

\usepackage{amssymb}
\usepackage{algorithm}
\usepackage{algorithmic}

\usepackage{graphicx}

\usepackage{xcolor,soul}

\usepackage{savesym}
\savesymbol{AND}

\DeclareMathOperator*{\mmz}{\mathrm{minimize}}
\DeclareMathOperator*{\sjt}{\text{subject to }}

\author{%
 \Name{Jia-Jie Zhu} \Email{jzhu@tuebingen.mpg.de}\\
 \addr Empirical Inference Department\\
 Max Planck Institute for Intelligent Systems, T\"ubingen, Germany
 \AND
 \Name{Moritz Diehl} \Email{moritz.diehl@imtek.uni-freiburg.de}\\
 \addr Department of Microsystems Engineering\\
 University of Freiburg, Freiburg, Germany%
 \AND
 \Name{Bernhard Sch\"olkopf} \Email{bs@tuebingen.mpg.de}\\
 \addr Empirical Inference Department\\
 Max Planck Institute for Intelligent Systems, T\"ubingen, Germany
}

\begin{document}

\maketitle

\begin{abstract}%
In this paper, we apply kernel mean embedding methods to sample-based stochastic optimization and control.
Specifically, we use the reduced-set expansion method as a way to discard sampled scenarios.
The effect of such constraint removal is improved optimality and decreased conservativeness.
This is achieved by solving a distributional-distance-regularized optimization problem. We demonstrated this optimization formulation is well-motivated in theory, computationally tractable, and effective in numerical algorithms. 
\end{abstract}


\section{Introduction}
\begin{figure}[tb!]
	\centering
	\vspace{-0.4cm}
	\includegraphics[width=0.375\columnwidth]{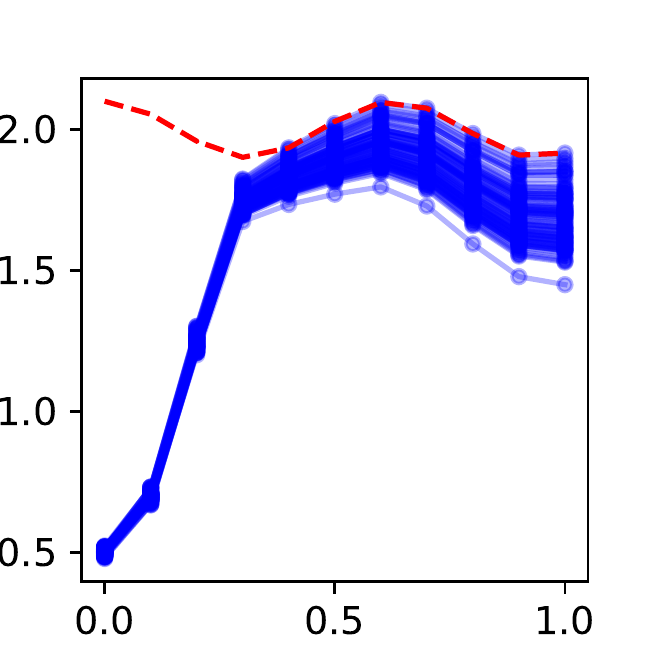}
	\includegraphics[width=0.375\columnwidth]{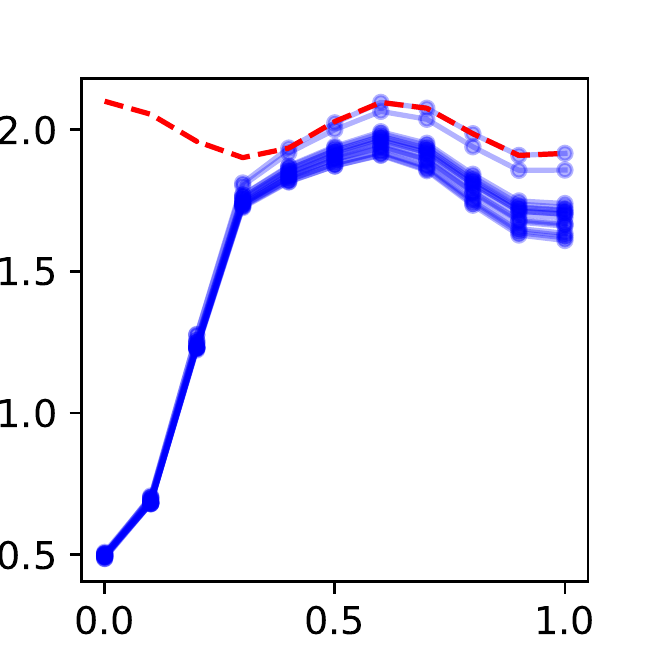}
	\put(-260,-10){\small{time step}} 
	\put(-340,60){\small{\rotatebox{90}{$x_1$}}}
	\caption{Steering the state trajectories (blue) without crossing constraints (red). Two figures depict optimization over different number of scenarios ($N=100 \text{ vs. } N=17$).}
	\vspace{-0.35cm}
	\label{fig:intro}
\end{figure}
Robustification against uncertain events is at the core of modern optimization and control. From the classic S-lemma to recent advances in distributionally robust optimization, we have witnessed computational tools giving rise to new robustification designs.
The classic ``worst-case'' approach sets out to robustify constraints against all realizations of disturbances in a mathematical model, resulting in the often over-conservativeness. 
Consider the illustrative example of a sample-based approach to solving constrained stochastic control problem in Figure~\ref{fig:intro}. In this case, we sample different numbers of realizations of the uncertainty in (left) and (right), and then solve the control problem under those realized scenarios.
Intuitively, one may expect the controller associated with more scenarios in (left) to be more robust against constraint violation than (right).
However, this results in conservative designs which may be reflected in high cost.

The central idea of modern data-driven robust optimization, in non-technical terms, is to use data samples to form empirical understanding of the true distribution. Then, one only seeks to robustify against this empirical understanding instead of the whole distribution support. One concrete relevance to our discussion, for example in Figure~\ref{fig:intro}, is that fewer scenarios translate to fewer constraints, which in turn lead to reduced cost. Therefore, this is a \emph{trade-off between optimality and feasibility}. 

In this paper, we show that constraint removal can be formulated as an optimization problem aiming to form a new distribution close to the empirical data distribution.
Our contributions are {\bf (1)} We formulate the constraint removal in stochastic programming and control as a tractable convex optimization problem with reproducing kernel Hilbert space (RKHS)-distance regularization or constraint. This formulation is well motivated in theory and effective in numerical studies. {\bf (2)} To our knowledge, this is the first use of RKHS-embedding reduced-set method in stochastic optimization and scenario approaches to control. 

\emph{Notation.}
In this work, symbol $\mathcal{H}$ denotes a reproducing kernel Hilbert space (RKHS). We write $\xi\sim P$ to denote that the random variable or vector (RV) $\xi$ follows the distribution law $P$.
By empirical distribution of the data, we mean the convex combination of Dirac-measures of the seen data $P_\mathrm{data}:=\frac1N\sum_i^N\delta(\xi_i)$ where $\{\xi_i\}_{i=1}^N$ is the data set.

\section{Background \& related work}
\label{sec:background}
\subsection{Stochastic programming and scenario optimization for control}
\label{sec:background_s}
In this paper, the problem of interest is the (chance-constrained) stochastic programming (SP; also known as stochastic optimization) in the following canonical formulation.
\begin{equation}
    \underset{x}{\mmz}\  \mathbb{E}[F(x, \xi)]\quad \sjt\  \operatorname{Pr}\{C(x, \xi) \leq 0\} \geq 1-\alpha.
    \label{eq:cc}
\end{equation}
As $\xi$ is assumed to be an RV, program~\eqref{eq:cc} may be intuitively understood as making decision $x$ under uncertainty originated from $\xi$ .
We consider the following sample-based SP (a.k.a. scenario approach).

Suppose we have a set of i.i.d. realizations $\{\xi_i\}_{i=1}^N$ of $\xi$, we solve the sample-based program
\begin{equation}
\underset{x}{\mmz}\  \mathbb{E}[F(x, \xi)]
\quad \sjt\  C(x, \xi_i) \leq 0\ \mathrm{for}\ i=1,\dots, N.
\label{eq:scenario}
\end{equation}
If $F$ and $C$ are convex in $x$, measurable in $\xi$, it can be shown that this formulation is a convex approximation to the original SP~\eqref{eq:cc}. 
As $N\to\infty$, the solution recovers that of the SP with level $\alpha=0$. 
However, with a large $N$, the solution to \eqref{eq:scenario} is overly conservative--- it aims to satisfy the constraints \emph{almost everywhere} in the distribution of $\xi$.
Therefore, the size $N$ trades off the conservativeness with constraint-satisfaction. Extensive research (e.g.,\cite{calafiore2006scenario,dentcheva2000concavity,luedtke2010integer} ) has focused on approaches to remove a subset of sampled constraints to reduce conservativeness of the solution.

Relevant to this paper, \cite{campi2011sampling} established guaranteed bounds for the constraint satisfaction probability and the number of removed constraints $\kappa$. Our method is built upon their sampling-and-discarding framework.
\cite{campi2013random} used $l1$ regularization to encourage sparsity in decision variables, which is different from our sparsity in RKHS expansion terms.
For readers who are interested in sample-based stochastic programming, good text references are given by Ch.5 of \cite{shapiro2009lectures} and Ch.9 of \cite{birge2011introduction}. 

\subsection{Reproducing kernel Hilbert space (RKHS) embeddings}\label{sec:bg_kernel}
A positive definite kernel is a real-valued bivariate, symmetric function $k(\cdot,\cdot):\mathcal{X}\times\mathcal{X}\to\mathbb{R}$ such that $\sum_{i,j=1}^n\alpha_i\alpha_jk(x_i,x_j) \geq 0$ for any $n\in\mathbb{N}$, $(\alpha_1,\ldots,\alpha_n)\in\mathbb{R}^n$, and $(x_1,\ldots,x_n)\in\mathcal{X}^n$. 
One may intuitively think $k(x,x')$ as a generalized similarity measure (inner product) between $x$ and $x'$ after mapping them into the feature space $\mathcal{H}$, $ k(x, x') = \langle \phi(x),\phi(x')\rangle_{\mathcal{H}}.$
We refer to \(\phi\) as \emph{feature map} associated with the kernel $k$, and $\mathcal H$ the associated RKHS.
A canonical kernel is the Gaussian kernel
$
	k(x,x') = \exp\left(-\frac{1}{2\sigma^2}\|x-x'\|_2^2\right)
$
where $\sigma > 0$ is a bandwidth parameter.

RKHS embedding, or \emph{kernel mean embedding} (KME) [\cite{Smola07Hilbert}] maps probability distributions to (deterministic) elements of a Hilbert space.
Mathematically, the KME of a random variable $X$ is given by the function
$\mu_X(\cdot) = \int  \langle\phi(x), \phi(\cdot)\rangle \ dP(x)$, which is a member of the RKHS.
For example, the RKHS associated with the second-order polynomial kernel consists of quadratic functions whose coefficients preserve statistical mean and variance. Gaussian kernel embeddings, on the other hand, preserve richer information up to infinite order.


\emph{Reduced-set expansion method using RKHS embeddings.}
Given a data set ${x_i}_{i=1}^N$, the sample-based KME is given by
$\hat \mu _ {x}   = \sum_{i=1}^N \alpha_i \phi (x_i),
$ where one can simply choose $\alpha_i=\frac1N$. It has been shown that one may use fewer than the total $N$ data sample to represent the distribution. This is the idea of reduced-set approximation. (cf. \cite{scholkopf2002learning}) Mathematically, this method seeks to find an embedding with fewer expansion terms
$\hat \mu^R_x =   \sum_{i=1}^{N_R} \alpha^i \phi (x_i)\approx \hat \mu _ {x}, N_R<N,$where the approximation is in the sense of RKHS distance measure. The reduced-set method forms the backbone of our approach.
We also note that there are other related approximation methods such as those of \cite{chen2012super, bach2012equivalence}.
Recently, \cite{zhu2019new} considered recursive applications of reduced-set method to uncertainty in stochastic systems. 

\section{Method}
\label{sec:method}
\subsection{Stochastic programming with reduced-set expansion of RKHS embeddings}
\label{sec:method_reduce}
We consider the sample-based formulation of the stochastic programming problem~\eqref{eq:scenario}.
Our main idea is to perform constraint removal systematically using the aforementioned RKHS embedding reduced-set methods. Typically, constraint removal discards low-probability scenarios to reduce conservativeness of the resulting solution.
Given a set of realized scenarios $\xi:=\{\xi_1, \dots \xi_n\}$ and positive definite kernel $k$, we formulate optimization problem as
\begin{equation}
\mmz _ \alpha  \| w^\top\alpha \|_1  \quad
     \sjt\ \| \sum_{i=1}^N \alpha_i \phi(\xi_i) - \hat \mu_{\xi}\|_{\mathcal H}\leq \epsilon,
\label{eq:l1}
\end{equation}
where $w\in\mathbb R^N$ denotes scaling vector for the $l1$-penalty. This can often be set to reflect specific concerns, such as the distance of states to the constraint.
The KME expansion weights $\alpha:=(\alpha_1,\dots,\alpha_N)\in \mathbb R^N$ need not sum to one.
$\hat\mu_\xi$ denotes the empirical KME estimator of the distribution $\hat \mu_{\xi}:=\frac1N\sum_{i=1}^N \phi(\xi_i)$.  
We further write down the equivalent Lagrangian form.
\begin{equation}
    \begin{aligned}&\mmz _ \alpha& & \| \sum_{i=1}^N \alpha_i \phi(\xi_i) - \hat \mu_{\xi}\|_{\mathcal H}^2 +\lambda\| w^\top\alpha \|_1.&
    \end{aligned}
    \label{eq:l1_reg}
\end{equation}
The resulting solution $\alpha^*$ is sparse due to the sparsity-inducing $l1$ term. We then discard the points~$\xi_i, i\in\mathcal I$ with the index set $\mathcal I = \{i\ |\ \alpha_i=0, i=1,\dots,n\}$. Finally, the we re-solve the stochastic programming problem with the reduced-set scenarios
\begin{equation}
\mathcal R := \{1,\dots,n\}\setminus\mathcal I.
\label{eq:index}
\end{equation}

The \emph{intuition} of the optimization formulation \eqref{eq:l1} and \eqref{eq:l1_reg} is to produce a subset of data whose distribution is close to the empirical data in the sense of RKHS-embedding distance $\| \hat \mu^R_{x} - \hat \mu_{x}\|_{\mathcal H}$. Meanwhile, the weighted $l1$-penalty incentivizes the solution to become sparse. Therefore, the solution to \eqref{eq:l1_reg} discards the ``corner" cases while maintaining the statistical information.
We outline the algorithmic procedure in Algorithm~\ref{alg:reduced}. 
\begin{algorithm}[tbp]
	\caption{RKHS approximation to stochastic programming}
	\label{alg:reduced}
	\begin{algorithmic}[1]
		\STATE Solve the sample-based stochastic programming problem~\eqref{eq:scenario}.
		\STATE Find the reduced-set RKHS embeddings
		\begin{equation*}
			\hat \mu^R_{x} =   \sum_{i\in\mathcal R} \alpha_i \phi(\xi_i).
		\end{equation*}
		by solving the convex optimization problem \eqref{eq:l1}. $\mathcal R$ is the reduced index set defined in \eqref{eq:index}. 
        \label{step:rkhs_reduce}
        \STATE Solve the stochastic programming problem according to the reduced set RKHS approximation (i.e., constraint removal by sparse optimization).
        \begin{equation}
			\mmz_x L(x), \quad \sjt   C(x, \xi_j) \leq 0,\ j\in \mathcal R.
		\end{equation}
		\vspace{-0.5cm}
		\label{step:redo_stoprog}
		\STATE {\bf Output:} Solution of the above reduced stochastic program.
	\end{algorithmic}
\end{algorithm}
The following lemma shows formulation~\eqref{eq:l1_reg} is computationally tractable.

If $\mathcal H$ in problem \eqref{eq:l1_reg} is the RKHS associated with a positive definite kernel, the objective of optimization problem is convex.
To see this, we use the sample based estimator for the RKHS distance to rewrite the objective as
\begin{equation}
    \begin{aligned}&\mmz _ \alpha& &  \alpha^\top K \alpha + \alpha^\top K\beta + \lambda\cdot \| w^\top\alpha \|_1.
    \end{aligned}
    \label{eq:l1_estimate}
\end{equation}
$\beta=(\frac1N,\dots,\frac1N)$ is a constant vector. 
$K:=\{k(x_i, x_j)\}_{i,j}$ is the gram matrix associated with the positive definite kernel, which implies $K\geq 0$. As $\|\cdot \|_1$ is convex, so is the optimization problem.

\begin{remark}    
We can equivalently write the constraint of program~\eqref{eq:l1} in the form of maximum mean discrepancy, $\| P^R - \hat P_\mathrm{data}\|_{\mathrm{MMD}}\leq \epsilon$.
Optimization problems with constraints on the probability distributions are often referred to as \emph{moment problems}.
\end{remark}
\subsection{Application to stochastic optimal control}
Let us  consider the following sample-based (scenario) formulation of stochastic optimal control problem (OCP).
\begin{equation}
    \begin{aligned}& \mmz_{u} & & \frac{1}{N}\sum_{i=1}^{N}\int_0^T L(x^i(t), u(t))\ dt\\
    & \text{subject to} & & \dot x^i = f(x^i, u, 
    \xi^i),\\& & & h(x^i, u,\xi^i) \leq 0,\\&&& x^i(0) = \xi_0^i, i=1,2,\dots, N.
    \end{aligned}
    \label{eq:ocp_sample}
\end{equation}
where $\xi, \xi_0$ are uncertain variables and $\xi^i, \xi_0^i$ their realizations.
The uncertainty in the initial state $\xi_0$ is particularly relevant to MPC designs.
After proper transcription and discretization, this OCP subsequently becomes the same form as the sample-based SP~\eqref{eq:scenario}, solvable by Algorithm~\ref{alg:reduced}.

\label{sec:theory}

\section{Numerical Experiments}

\subsection{Min-max robust regression}
We first consider a synthetic stochastic programming problem given in the form of the following min-max robust regression.
\begin{equation}
\begin{array}{llll}
	\underset{x}{\ \ \mmz}\ S, & {{\sjt }}{\left|A_{i} x-b_i\right| \leqslant S}, \forall i.
\end{array}
\label{eq:minmax}
\end{equation}
For simplicity, we consider scalars $A_{i}$ and $b_i$ generated randomly according to the distributions.
$
n_1,n_2\sim\mathrm N(0, 1),\,
A_i = 3 + 3n_1,\,
\ b_i = A_i x^* + 5n_2.
$
where $x^*$ is the (unknown) true parameter drawn from Uniform($[2,3]$).
\begin{figure}[tb!]
	\centering
	\includegraphics[width=0.9\columnwidth]{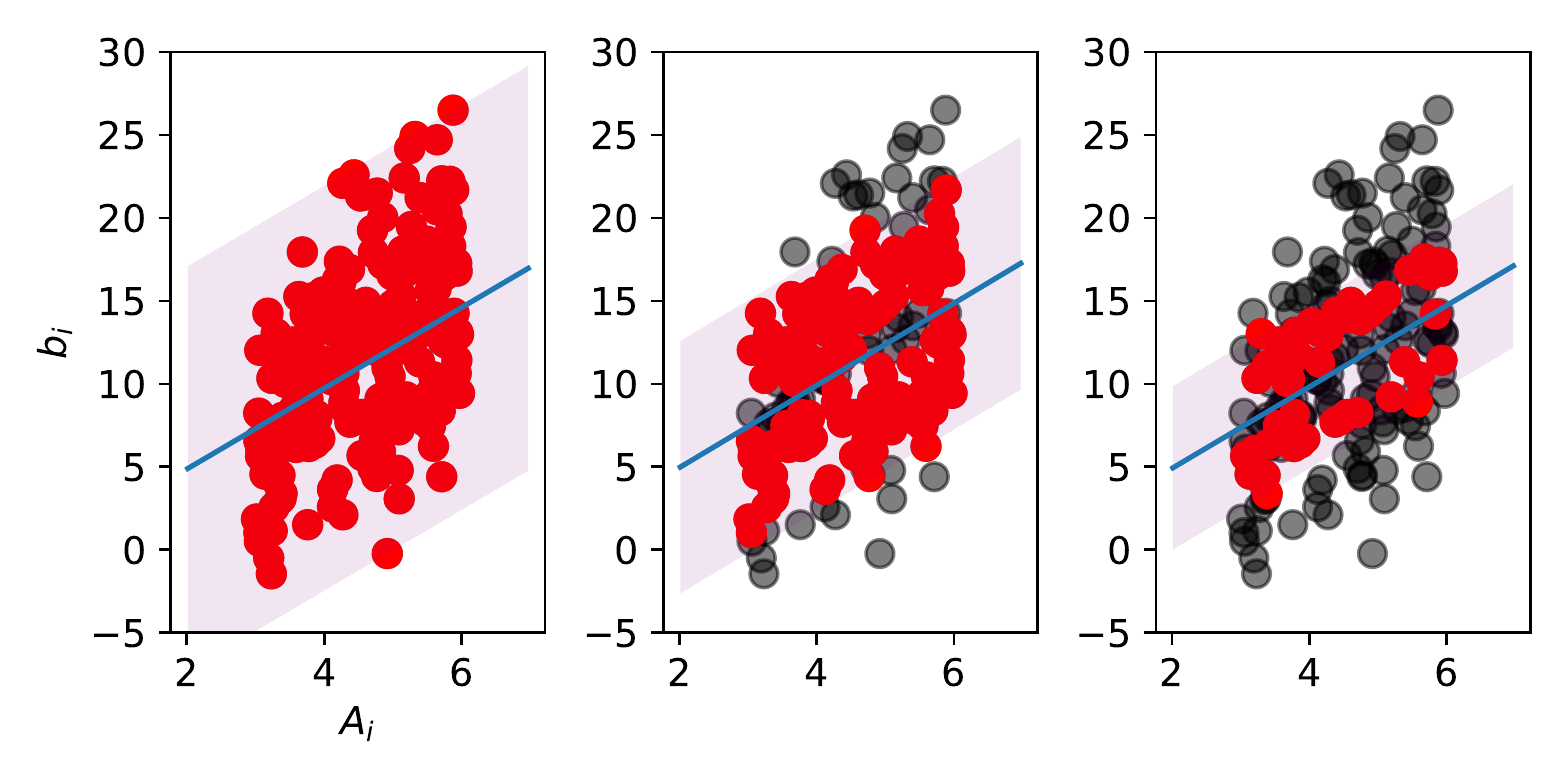}
	\vspace{-0.3cm}
	\caption{Solutions of the min-max robust regression problem. Three figures correspond to three different regularization coefficients $\lambda$ and number of discarded scenarios $\kappa$: (left) $\lambda=0,\kappa=0$. (center) $\lambda=0.01, \kappa=57$. (right) $\lambda=0.05, \kappa=144$.
	The shaded strip denotes the robust margin. Red points are the selected points by Step~\ref{step:rkhs_reduce} in Algorithm~\ref{alg:reduced}. Dark points correspond to discarded scenarios and constraints. We used Gaussian kernel of bandwidth $\frac{1}{\sqrt{2}}$ to calculate the RKHS embedding in Algorithm~\ref{alg:reduced}.
	}
	\vspace{-0.7cm}
	\label{fig:lasso3}
\end{figure}

\begin{figure}[tb!]
	\centering
	\includegraphics[width=0.75\columnwidth]{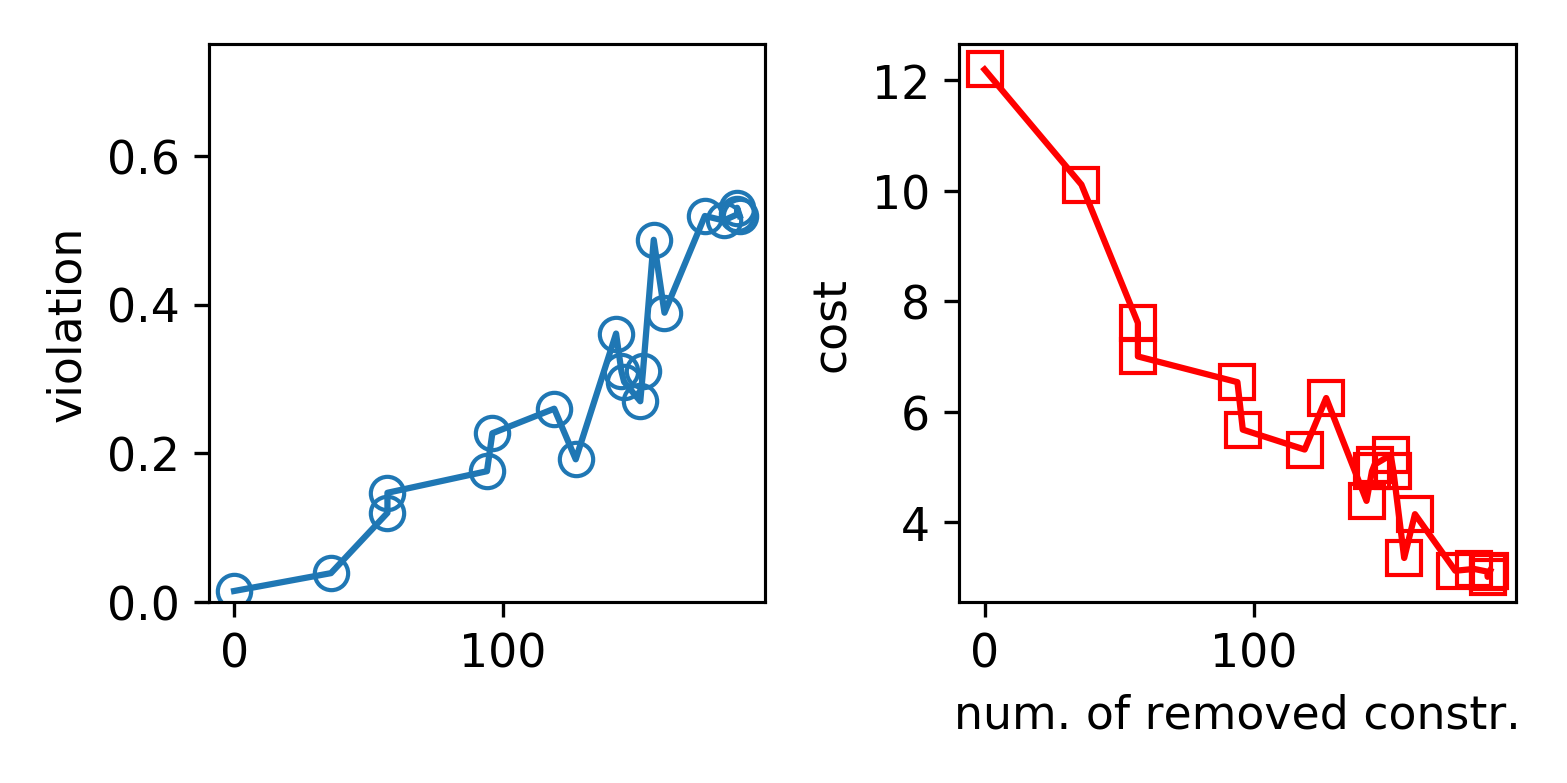}
	\vspace{-0.3cm}
	\caption{Estimates produced by a large-sample ($N=1000$) Monte Carlo simulation.
	(left). Constraint violation probability of the solution produced by Algorithm~\ref{alg:reduced}. This is estimated by the Monte Carlo estimation $\frac1N\sum_{i=1}^{N_R} \mathbb I (\left|A_{i} x-b_i\right| >  S)$ where $I$ is the indicator function of random events.
	(right).The new expected cost associated with the solution produced by Algorithm~\ref{alg:reduced}. This is estimated by the Monte Carlo estimation $\frac1N\sum_{i=1}^{N_R} F(x^*, \xi_i)$ where $x^*$ is the solution by the proposed method.
	}
	\vspace{-1cm}
	\label{fig:cost_vs_vio}
\end{figure}

Given the computed solution to the full program~\eqref{eq:minmax} $\hat x$, let us consider the quantity of interest $\hat \xi_i:=A_i\hat x-b_i$, which is an RV due to the uncertainty in $A_{i}$ and $b_i$. We now apply Algorithm~\ref{alg:reduced} to find the reduced-set embedding of $\{\hat \xi_i\}_i^{N_R}$, 
$
\hat \mu _ {\hat\xi}   = \sum_{i=1}^{N_R} \alpha^R_i \phi({\xi_i)}.
$ 
In step~\ref{step:redo_stoprog}, in solving program~\eqref{eq:l1_reg}, we used the scaling factor $w_i\propto\exp{(T\hat\xi_i)}$ to incentivize the removal of ``corner'' points ($T$ may be thought of as the ``softness'' parameter of this softmax scaling factor). We then remove the constraints with identified index set $\mathcal I$,
$
|A_i x- b_i| \leq S,\ \forall i\in\mathcal I,
$
from the stochastic program and re-compute a solution.
Following our discussion in the previous sections, this embedding captures the distribution information while discarding the rare scenarios. This is done by solving the sparse optimization problem~\eqref{eq:l1}. The results are illustrated in Figure~\ref{fig:lasso3}. As we can see, scenarios associated with ``corner" data points are not selected, causing the reduction in conservativeness. 
Figure~\ref{fig:cost_vs_vio} illustrates the effect of constraint removal on the optimal objective value and constraint violation. See the caption for detailed description.

\subsection{Stochastic control}
We now consider the Van der Pol oscilator model
\begin{equation}
	\frac{d}{dt}{\left[\begin{array}{l}   {x_{1}} \\  {x_{2}} \end{array}\right]}=\left[\begin{array}{c} {x_{2}} \\  {-0.1\left(1-x_{1}^{2}\right) x_{2}-x_{1}+u}    \end{array}\right].
	\label{eq:vanderpol}
\end{equation}
The goal of the control design is to steer the system state $x_1$ to a certain level. This is formulated as the following OCP.
\begin{equation}
	\begin{array}{ll}
	{\underset{x(\cdot), u(\cdot)}{\operatorname{minimize}}} & {\int_{0}^{T}\left\|x_{1}(t)-3\right\|_{2}^{2} d t} \\
	{\text { subject to }} & {\dot{x}(t)=f(x(t), u(t)) \quad \forall t \in[0, T]} \\
	{} & {-40 \leq u(t) \leq 40 \quad \forall t \in[0, T]} \\
	{} & {-0.25 \leq x_{1}(t) \leq 2+0.1 \cos (10 t) \quad \forall t \in[0, T]} \\
	{} & {x(0)=s}
	\end{array}
	\label{eq:ocp2}
\end{equation}

We sample the i.i.d. uncertainty realizations $\{s_1, \dots s_n\}\sim N(m,\Sigma)$, where $m=[0.5\ \ 0]^{T},\Sigma=\left[\begin{array}{cc}
	{0.01^{2}} & {0} \\
	{0} & {0.1^{2}}
	\end{array}\right]$.
Because of the nonlinear dynamics, we cannot propagate the uncertainty in a tractable manner as in LQG without resorting to approximations. We use the sampled scenarios to form the OCP~\eqref{eq:ocp_sample}. The continuous-time dynamics is transcribed using multiple-shooting with CVodes (interfaced with CasADi) integrator. We then solve the discretized OCP with IPOPT to obtain the optimal control. An example of states associated with the solution is given in Figure~\ref{fig:intro} (left). The total time horizon is $1.0$s and we consider $10$ control steps in this experiment.

Let us consider the quantity of interest $ \xi_i(t):=2+0.1 \cos (10 t) - x^i_{1}(t)$, the distance from state position to the upper bound constraint. This quantity reflects how close we are to be infeasible. It is random due to the states being a function of RV $s$. In Step~\ref{step:rkhs_reduce} of Algorithm~\ref{alg:reduced}, the scaling factor is taken to be $w_i\propto\exp{(\frac{C}{\epsilon+\min_t(\hat\xi_i(t))})}$ to encourage the removal of close-to-constraint trajectories.

We are now ready to apply Algorithm~\ref{alg:reduced} to find the reduced-set embedding of $\{\hat \xi_i(t)\}_i^{N_R}, t=1,\dots T$, 
$
\hat \mu _ {\xi}   = \sum_{i=1}^{N_R} \alpha_i \phi({\xi_i)}
$, where $\xi_i$ is a vector comprising $\xi_i(t)$ at all time steps.
Finally, we re-solve the subsequent reduced-set SP---OCP.
Figure~\ref{fig:intro} (right) illustrates the reduced number of scenarios.

\begin{figure}[b!]
	\vspace{-0.3cm}
	\centering
	\includegraphics[width=0.375\columnwidth]{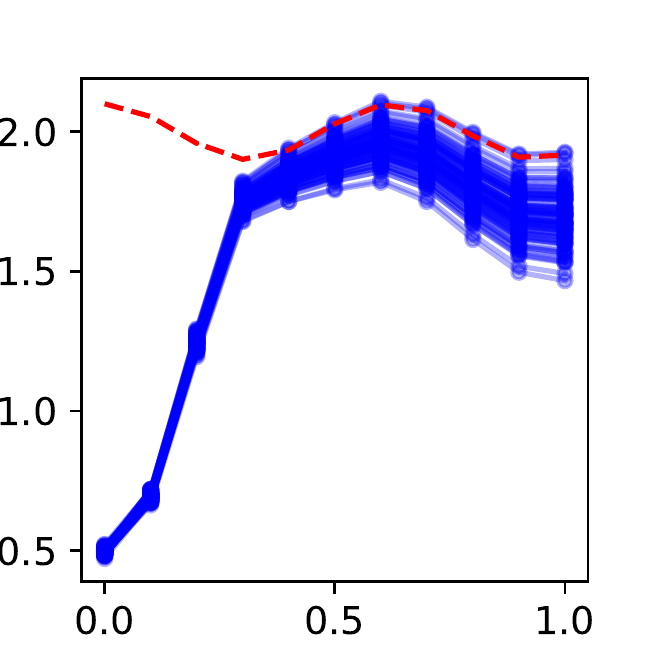}
	\includegraphics[width=0.375\columnwidth]{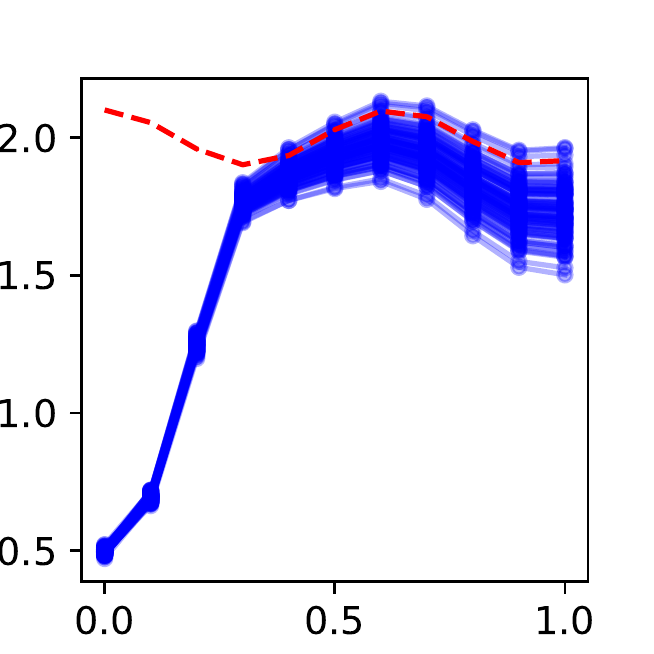}
	\put(-260,-10){\small{time step}} 
	\put(-340,60){\small{\rotatebox{90}{$x_1$}}}
	 
	\caption{(left) The optimistic controller evaluated in independent Monte Carlo simulations with $N=100$ trajectories. The controller is produced with Algorithm~\ref{alg:reduced} with $\kappa=14$ removed scenarios ($l1$ regularization coefficient $10^{-5}$). The estimated constraint-violation probability
	is $2\%$.
	(right) Controller produced with $\kappa=83$ removed scenarios ($l1$ regularization $5\times 10^{-3}$). The constraint violation probability is $5\%$.
	}
	\vspace{-0.5cm}
	\label{fig:lasso_mpc}
\end{figure}

\begin{figure}[hbt!]
	\centering
	\includegraphics[width=0.75\columnwidth]{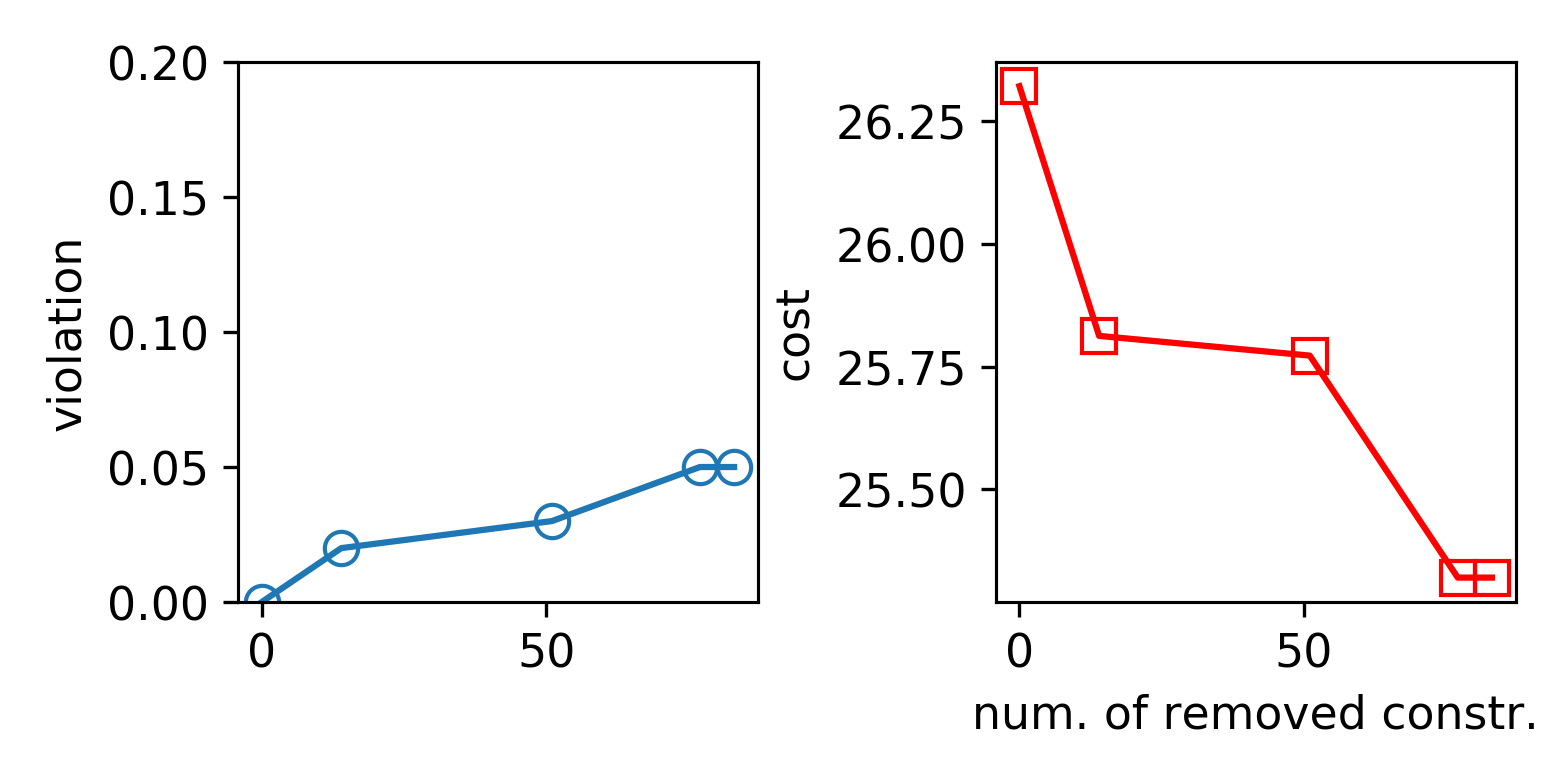}
	\caption{Performance estimates produced by independent ($N=100$) Monte Carlo simulations.
	(left). Constraint violation probability of the solution produced by Algorithm~\ref{alg:reduced}. This is estimated by the Monte Carlo estimation $\frac1N\sum_{i=1}^{N_R} \mathbb I (C(x,u) > 0)$ where $I$ is the indicator function. $C(x,u)\leq 0 $ denotes all the constraints in OCP~\eqref{eq:ocp2}.
	(right).The new expected cost of OCP associated with the solution produced by Algorithm~\ref{alg:reduced}. This is estimated by the Monte Carlo estimation $\frac1N\sum_{i=1}^{N} \int_{0}^{T} L(x^i(t), u^*(t))\ dt$ where $u^*$ is the optimistic controller.
	}
	\label{fig:mpc_cost_vs_vio}
\end{figure}

After we applied Algorithm~\ref{alg:reduced}, we obtain the ``optimistic'' controller $u^*$. 
To evaluate the performance of this controller, we use a large-sample Monte Carlo simulation to estimate the constraint violation probability, i.e., $\operatorname{Pr}\{C(x, \xi) \leq 0\}$ in the chance-constrained SP~\eqref{eq:cc}, as well as the expected cost over the large-sample simulation $\frac1N\sum_{i=1}^{N} \int_{0}^{T} L(x^i(t), u^*(t)) dt$. 
We plot the state trajectories with different number of removed constraints in Figure~\ref{fig:lasso_mpc}. The trade off between those is illustrated in Figure~\ref{fig:mpc_cost_vs_vio}.
The result makes intuitive sense that the more constraints we remove, the less the conservativeness, but with higher violation probability. 
 

\section{Discussion}
This paper proposed a distributional-distance-regularized optimization formulation for stochastic programming under the framework of sampling-and-discarding. We demonstrated effective conservativeness reduction in data-driven optimization and control tasks. Although we did not study the guaranteed bounds, all analysis in \cite{campi2011sampling} applies to our case.
Our on-going work is to apply our approach to distributionally robust optimization and control design.

\acks{We thank Joris Gillis for his help in the numerical experiment, Wittawat Jitkrittum for his helpful feedback. This project has received funding from the European Union’s Horizon 2020 research and innovation programme under the Marie Skłodowska-Curie grant agreement No 798321, the German Federal Ministry for Economic Affairs and Energy (BMWi) via eco4wind (0324125B) and DyConPV (0324166B), and by DFG via Research Unit FOR 2401.
}

\bibliography{l4dc2020} 

\end{document}